\documentclass[12pt]{amsart}
\usepackage{amsmath,amssymb,amsfonts,mathrsfs,exscale,colonequals,a4wide,color,mathtools,url}

\usepackage[curve,matrix,arrow,cmtip]{xy}
\usepackage{verbatim}

\usepackage{cite}

\usepackage[utf8]{inputenc} 

\newdir^{ (}{{}*!/-3pt/\dir^{(}}   
\newdir_{ (}{{}*!/-3pt/\dir_{(}}   
\entrymodifiers={+!!<0pt,\fontdimen22\textfont2>}

\let\phi\varphi
\let\epsilon\varepsilon

\let\leq\leqslant
\let\geq\geqslant

\newtheorem{Thm}{Theorem}[section]

\newtheorem{Lem}[Thm]{Lemma}

\numberwithin{equation}{section}

\def\qed{{\hskip0pt\unskip\unskip\nobreak\hfil\penalty50
          \hskip1em\hbox{}\nobreak\hfil
           {$\square$}
          \parfillskip=0pt\finalhyphendemerits=0
          \par}\medskip}

\newenvironment{Proof}
               {\noindent{\bf Proof.}\ }
               {\qed}

\newenvironment{Proofof}[1]
               {\noindent{\bf Proof of #1.}\ }
               {\qed}

\newcommand{\BN}{{\mathbb{N}}}

\newcommand{\BZ}{{\mathbb{Z}}}

\newcommand{\CB}{{\mathcal B}}

\newcommand{\CF}{{\mathcal F}}

\begin{document}

\title{Normality in Pisot Numeration Systems}

\author{%
Adrian-Maria Scheerer%
}
\thanks{TU Graz, \tt scheerer@math.tugraz.at}

\begin{abstract}
Copeland and Erd\H{o}s \cite{copeland1946} showed that the concatenation of primes when written in base $10$ yields a real number that is normal to base $10$. We generalize this result to Pisot number bases in which all integers have finite expansion.
\end{abstract}

\date{\today}

\maketitle


%
%

\section{Introduction}
Let $x$ be a real number and $b \geq 2$ a positive integer. Then $x$ has a $b$-adic representation of the form
\begin{equation*}
x = \lfloor x \rfloor +  \sum_{i=1}^\infty \epsilon_i b^{-i}
\end{equation*}
where $\epsilon_i \in \{0,1, \ldots b-1\}$ are the digits of $x$ and $\lfloor x \rfloor$ denotes the integer part. We call $x$ \emph{normal to base $b$}, if any block $d=d_1 d_2 \ldots d_k$ of $k\geq1$ digits occurs with the expected frequency in the $b$-adic representation of $x$. This means that
\begin{equation*}
\lim_{n\rightarrow \infty} \frac{1}{n} N_d(x,n) = \frac{1}{b^k},
\end{equation*}
where $N_d(x,n)$ counts the occurrences of the block $d$ within the first $n$ digits of $x$. 
A real number $x$ is called \emph{absolutely normal} if it is normal to every base $b\geq2$.\\

The terminology of a normal number can be extended to the context when the underlying base is no longer an integer. R\'{e}nyi \cite{Renyi} introduced and Parry \cite{Parry} studied numeration systems with respect to real bases $\beta > 1$. Each real number $x$ has a representation of the form
\begin{equation*}
x = \sum_{i=L}^{-\infty} \epsilon_i \beta^i,
\end{equation*}
with digits $\epsilon_i \in \{0,1, \ldots, \lceil \beta \rceil -1\}$. One way to produce the digits is the so-called \emph{greedy algorithm} using the transformation $T_\beta: x \mapsto \beta x \pmod 1$ on the unit-interval. In a natural way, $T_\beta$ corresponds to the shift-operator on the set $W^\infty$ of right-infinite sequences over $\{0,1, \ldots, \lceil \beta \rceil -1\}$ and each $x$ corresponds to its sequence of digits. A sequence $\omega$ in $W^\infty$ is called \emph{$\mu$-normal} for a given shift-invariant measure $\mu$ on $W^\infty$, if all possible finite patterns of digits occur in $\omega$ with asymptotic frequency given by $\mu$. Consequently, the real $x$ is called \emph{$\mu$-normal} if its sequence of digits is $\mu$-normal. Details will be made clear in the next section.\\


From a modern approach, using Birkhoff's point-wise ergodic theorem, it is immediate that almost all numbers are normal to a fixed \footnote{Borel's result also follows, but for the moment we want to restrict the discussion to one single base.} base $b$. The map $ T_{b} : x \mapsto b x \pmod 1$ on the unit-interval is the underlying ergodic transformation which preserves Lebesgue-measure. Knowing this, the existence of normal numbers to base $b$ is in a certain sense not very surprising. However, in this context, there are two observations that make the study of normal numbers interesting.\\


\emph{First:} The explicit construction of normal numbers. The study of normal numbers dates back to Borel \cite{Borel}, who in 1909 showed that Lebesgue-almost all numbers are absolutely normal. However, the first explicit example of a normal number is due to Champernowne \cite{Champernowne} in 1933. He showed that the concatenation of integers, when written in base $10$,
\begin{equation*}
0,1\ 2\ 3\ 4\ 5\ 6\ 7\ 8\ 9\ 10\ 11\ 12\ 13 \ldots,
\end{equation*}
is normal to base $10$. Copeland and Erd\H{o}s \cite{copeland1946} showed that
\begin{equation*}
0, 2\ 3\ 5\ 7\ 11\ 13\ 17\ 19 \ldots,
\end{equation*}
i.e. the concatenation of primes in base $10$ is normal to base $10$. Besicovitch \cite{Besicovitch} showed that the decimal formed by concatenation of the squares in base $10$ is normal to base $10$. This construction was extended to general integer-valued polynomials by Davenport and Erd\H{o}s \cite{DavenportErdos} and by Schiffer \cite{Schiffer} and Nakai and Shiokawa \cite{MR1064444}, \cite{MR1197421} to more general polynomial settings. Nakai and Shiokawa \cite{MR1472814} also evaluated polynomials at primes, and Madritsch \cite{MR3273499} showed that numbers generated by pseudo-polynomial sequences along the primes are normal. Further constructions of normal numbers in the spirit of Copeland and Erd\H{o}s and Erd\H{o}s and Davenport include \cite{MR3032395} and \cite{MR2406483}.

These constructions, most notably the one due to Champernowne, have subsequently been generalized to other number systems. To mention are the works by Ito and Shiokawa \cite{ito1975} who generalized the Champernowne-construction to real bases $\beta > 1$, and by Madritsch and Mance \cite{MadritschMance} who modified the construction to produce normal sequences in general symbolic dynamical systems. Bertrand-Mathis and Volkmann \cite{BertrandMathisVolkmann} give a generalized Copeland-Erd\H{o}s construction to symbolic dynamical systems.\\

In this paper we prove a polynomial Copeland-Erd\H{o}s-construction to bases which are not integers. We use results from the work of Bertrand-Mathis and Volkmann \cite{BertrandMathisVolkmann} which in turn is extending the original work of Copeland and Erd\H{o}s \cite{copeland1946}.\\

\emph{Second:} Not \emph{all} numbers are normal (almost all cannot be improved to all, hence $T_{b	}$ is not uniquely ergodic), so the ergodic theorem is strict. This can be seen in context of the following `test' for ergodicity.

\begin{Thm}[See Theorem 1.4 in \cite{BillingsleyErgodicTheory}] \label{Test}
Let $\CF_0$ be a field\footnote{The underlying sigma-field $\CF$ is in our case the Borel sigma-algebra which is generated by $b$-adic intervals.} generating $\CF$. If $T$ respects\footnote{Almost every orbit under $T$ visits $A$ with the expected asymptotic frequency.} every $A$ in $\CF_0$, then $T$ is ergodic.
\end{Thm}

Hence normal numbers, or rather the existence of non-normal numbers, can be used to test the underlying transformation for unique ergodicity. Although the Lebesgue-measure is not the only $T_{b}$-invariant probability measure on the unit-interval, it is the only one that maximizes entropy (see the following section). Accepting this for the moment as a definition of uniqueness, it is possible to study normality in greater generality with respect to a given transformation that can be different from $T_b$.

%
%


\section{Preliminaries}
In the following we present a condensed introduction to $\beta$-expansions, Pisot numbers and symbolic dynamical systems - the context in which we want to state our result.\\

Let $\beta > 1$ be a fixed real number. A \emph{$\beta$-expansion} of a non-negative real number $x$ is a representation of $x$ as a sum of integer powers of $\beta$ of the form
\begin{equation}\label{expansion}
x = \sum_{i=L}^{-\infty} \epsilon_i \beta^i,
\end{equation}
where the digits $\epsilon_i \in \{0,1, \ldots \lceil \beta \rceil -1\}$ are obtained by the following \emph{greedy algorithm}. Let $L\in \BZ$ such that $\beta^L \leq x < \beta^{L+1}$ and put $\epsilon_L = \lfloor x / \beta^L \rfloor$ and $r_L = \{x/ \beta^L \}$. For $L \geq i > - \infty$, define recursively $\epsilon_i = \lfloor \beta r_{i+1} \rfloor$ and $r_i = \{\beta r_{i+1} \}$. $\beta$-expansions have been introduced and studied by R\'{e}nyi \cite{Renyi} and Parry \cite{Parry}. 

Let $T_\beta$ be the \emph{$\beta$-transformation} $T_\beta : [0,1) \rightarrow [0,1)$, $x \mapsto \{ \beta x \}$. The digits in the $\beta$-expansion \ref{expansion} are given by $\epsilon_i = \lfloor \beta T_\beta^{i-1}(x) \rfloor$. R\'{e}nyi \cite{Renyi} showed that there is a unique normalized measure $\mu_\beta$ on $[0,1)$ that is invariant under $T_\beta$ and equivalent to the Lebesgue measure. This measure also maximizes the entropy of the corresponding symbolic dynamical system and we use it to define normal numbers in base $\beta$, see below.\\

A \emph{Pisot number} $\beta$ is a real algebraic integer $\beta > 1$ such that all its conjugates have absolute value less than $1$. For a Pisot number of degree $d$, we denote its conjugates by $\beta_i$, $i = 2, \ldots, d$, and the corresponding conjugations by $\sigma_i$.

We work with Pisot numbers such that every positive integer has finite $\beta$-expansion. A criterion for when this is the case can for example be found in \cite{frougny1992finite}.\\

In our notation concerning symbolic dynamical systems we follow Bertrand-Mathis and Volkmann \cite{BertrandMathisVolkmann}.

Let $A$ be a finite alphabet. 
$A^\ast$ is the set of all finite (possibly empty) words over $A$ and
$A^\BN$ is the set of all (right-)infinite words over $A$.
We call a subset $L$ of $A^\ast$ a language. 
$L^\ast$ denotes the set of all finite concatenations of words from $L$. 
Let $W(L^\ast)$ be the set of all non-empty factors of words in $L^\ast$. 
For a finite word $\omega$ we let $\Vert \omega \Vert$ be its length.
A language $L$ is said to be connecting of order $j \geq 0$ if for any two words $a,b \in W(L^\ast)$ there is a word $u = u(a,b) \in W(L^\ast)$ of length $j$ such that $aub \in W(L^\ast)$. For each $a, b  \in W(L^\ast)$ we choose one $u = u(a,b)$ and introduce the notation $a \oplus b := aub$. In the applications we have in mind, this intermediary word will simply consist of $0$'s.
We write $W(L^\ast) = \bigcup_{n\geq 1} L_n$ where $L_n$ is the subset of words in $W(L^\ast)$ of length $n$. Also denote by $L_n'$ the subset of $W(L^\ast)$ of words of length less or equal to $n$. 
For a language $L$ we denote by $W^\infty = W^\infty(L)$ the set of all infinite words generated by $L$, i.e. the set of all $\omega = a_1 a_2 \ldots$ such that $a_i a_{i+1} \ldots a_k \in L$ for all $1\leq i < k < \infty$.\\

We introduce the discrete topology on the alphabet $A$ and the corresponding product topology on the set of sequences $A^\BN$. With each language $L$ we associate the symbolic dynamical system
\begin{equation*}
S_L = (W^\infty, \CB, T, I),
\end{equation*}
where $W^\infty = W^\infty(L)$; $\CB$ is the $\sigma$-algebra generated by all cylinder sets of $A^\BN$, i.e. sets of the form
\begin{equation*}
c(\omega) = \lbrace a_1 a_2 \ldots \in A^\BN \mid a_1 a_2 \ldots a_n = \omega \rbrace
\end{equation*}
for some word $\omega \in A^\ast$ of length $n$. $T$ is the shift operator and $I$ is the set of all $T$-invariant probability measures $\mu$ on $\CB$. We will write $\mu(\omega)$ instead of $\mu(c(\omega))$ for a finite word $\omega$.\\

With each symbolic dynamical system $S_L$ we associate the entropy
\begin{equation*}
h(W^\infty) = \sup_{\mu \in I} h(\mu),
\end{equation*}
where $h(\mu)$ denotes the entropy of the measure $\mu$ (cf. Chapter 2 of Billingsley \cite{BillingsleyErgodicTheory}). For finite alphabets it is known that there always exists a unique measure $\mu$ that maximizes the entropy, i.e. such that $h(W^\infty) = h(\mu)$ (e.g. Proposition 19.13 in \cite{denker1976ergodic}). In the context of the symbolic dynamical system generated by a $\beta$-expansion the measure $\mu_\beta$ is precisely the measure with maximum entropy equal to $\log \beta$ (cf. \cite{Hofbauer1978} and the work by Bertrand-Mathis and Volkmann \cite{BertrandMathisVolkmann}). In the following, we will work with this maximal measure.\\

For an infinite word $\omega = \omega_1 \omega_2 \ldots \in W^\infty$ and a block $d=d_1d_2 \ldots d_k \in L$ we denote with $N_d(\omega,n)$ the number of occurrences of $d$ within the first $n$ letters of $\omega$. If the word $\omega$ is finite, we denote by $N_d(\omega)$ the occurrences of $d$ in it. An infinite word $\omega \in W^\infty$ is called \emph{$\mu$-normal} or \emph{$\mu$-normal sequence} if for all $d\in W(L^\ast)$
\begin{equation*}
\lim_{n\rightarrow \infty} \frac{1}{n} N_d(\omega,n) = \mu(d).
\end{equation*}
We note that a non-negative real number $x$ is \emph{$\mu_\beta$-normal} if and only if the sequence of digits in its $\beta$-expansion is a $\mu_\beta$-normal sequence.\\

In this terminology, the main result of Bertrand-Mathis and Volkmann \cite{BertrandMathisVolkmann} is the following

\begin{Thm}\label{MainThmABM}
Let $L$ be a connecting language and $a_1, a_2, \ldots$ a sequence of different elements of $W(L^\ast)$, $\Vert a_1 \Vert \leq \Vert a_1 \Vert \leq \ldots$, satisfying the generalised Copeland-Erd\H{o}s condition:
\begin{equation*}
\forall \epsilon > 0 \ \exists n_0(\epsilon) \ \forall n \geq n_0 \  \sharp \lbrace a_\nu \mid \Vert a_\nu \Vert \leq n \rbrace > \vert L_n'\vert^{1-\epsilon}.
\end{equation*}
Then the infinite word $a=a_1 a_2 \ldots \in W^\infty$ is normal.
\end{Thm}

The symbolic dynamical system \emph{generated by the $\beta$-shift} (or \emph{corresponding to the $\beta$-expansion}) arises naturally when viewing the $\beta$-expansions of real numbers $x\in[0,1)$ as infinite words over the alphabet $\{0,1,\ldots \lceil \beta \rceil\}$. $W^\infty$ is the set of these right-infinite sequences, $\CB$ the $\sigma$-algebra generated by all cylinder sets, $T$ the shift operator (it corresponds to the $\beta$-transformation $T_\beta$), and $I$ the set of all $T$-invariant probability measures on $\CB$. We work with the unique entropy-maximizing measure $\mu$ in $I$. It corresponds to $\mu_\beta$ on $[0,1)$ in the sense that for a finite word $\omega$ the measure of the cylinder set $\mu(\omega)$ is the same as the measure $\mu_\beta(\tilde{\omega})$ of the set $\tilde{\omega}$ of all real numbers in $[0,1)$ whose $\beta$-expansion starts with $\omega$. We allow us to speak of these two concepts interchangeably.\\

For a given Pisot number $\beta$, denote by $(n)_\beta$ the word over $\{0,1, \ldots \lceil \beta \rceil -1 \}$ that corresponds to the $\beta$-expansion of the positive integer $n$. We prove the following polynomial generalization of \cite{copeland1946}.

\begin{Thm}\label{ourThm}
Let $\beta$ be a Pisot number such that all integers have finite $\beta$-expansion and let the measure $\mu_\beta$ be as before. Let $f$ be a polynomial of degree $g$ that maps positive integers to positive integers. Then
\begin{equation*}
(f(2))_\beta \oplus (f(3))_\beta \oplus (f(5))_\beta \oplus (f(7))_\beta \oplus (f(11))_\beta \oplus \ldots
\end{equation*}
is a $\mu_\beta$-normal sequence. 
\end{Thm}

In the context of the dynamical system generated by the $\beta$-shift, the entropy is $\log \beta$. Hence by Lemma 2 of \cite{BertrandMathisVolkmann}, we have bounds on the number of words of length $n$ in $W(L^\ast)$. For $n$ sufficiently large we have $\beta^n \ll \vert L_n \vert \ll \beta^n$, where the implied constants do not depend on $n$. Therefore
\begin{equation}
\vert L_n' \vert = \sum_{\nu=1}^n \vert L_\nu \vert \ll \beta^n
\end{equation}
for all large $n$.

\section{General Case}

First we need upper and lower bounds of the length of the $\beta$-expansion of integers. Under the assumption of $\beta$ being a Pisot number such that all integers have finite $\beta$-expansions we can in fact show that the lengths of these expansions are asymptotically of logarithmic order of magnitude.

Note that if $n = \sum_{i=L(n)}^{-R(n)} \epsilon_i \beta^{i}$ we call $\sum_{i=0}^{L(n)} \epsilon_i \beta^{i}$ its integer part and $\sum_{i=-1}^{-R(n)} \epsilon_i \beta^{i}$ its fractional part. In the following we will think of $n$ as fixed and omit writing the dependency on it in the lengths $L$ and $R$.

\begin{Lem}\label{boundsgeneralcase}
Let $\beta$ be a Pisot number of degree $d$ such that all natural numbers have finite $\beta$-expansion. For the length $R(n)$ of the fractional part of $n$ upper and lower bounds of the following form hold, for sufficiently large $n$:
\begin{align*}
\delta \log n \leq R(n) \leq \delta' \log n
\end{align*}
where $\delta$ is a positive constant (specified in the proof) and the difference $\delta' - \delta>0$ can be chosen to be arbitrarily small.
\end{Lem}

\begin{Proof}
We have $\frac{n}{\beta^{L+1}} \in [0,1)$. 
Following an argument of Proposition 3.5 Frougny and Steiner \cite{FrougnySteiner}, for a certain number $k$ the number $T_\beta^k(\frac{n}{\beta^{L+1}})$ is an element of the finite set
\begin{equation*}
Y = \lbrace y \in \mathbb{Z}[\beta] \cap [0,1) \mid \vert \sigma_j(y) \vert < 1 + \frac{\lfloor \beta \rfloor}{1- \vert \sigma_j(\beta)\vert} \quad \text{for} \quad 2 \leq j \leq d \rbrace.
\end{equation*}
To see this, let the $\beta$-expansion of $n$ for the moment be $\epsilon_1 \epsilon_2 \ldots$ and put $z := \frac{n}{\beta^{L+1}}$. Then for all $k \geq 0$,
\begin{equation*}
T_\beta^k(z) = \beta T_\beta^{k-1}(z) - \epsilon_k = \ldots = \beta^k z - \sum_{l=1}^k \epsilon_l \beta^{k-l}.
\end{equation*}
Hence for all $k \geq 0$ and $2 \leq j \leq d$,
\begin{equation*}
\vert \sigma_j(T_\beta^k(z)) \vert = \vert \sigma_j(\beta)^k \sigma_j(z) - \sum_{l=1}^k \epsilon_l \sigma_j(\beta)^{k-l} \vert 
< \vert \sigma_j (\beta)\vert ^k \vert \sigma_j(z)\vert + \frac{\lfloor \beta \rfloor}{1 - \vert \sigma_j(\beta) \vert}.
\end{equation*}
Let $k$ be equal to
\begin{equation*}
\max_{2\leq j \leq d} \left\lceil - \frac{\log \vert \sigma_j(\frac{n}{\beta^{L+1}})\vert}{\log \vert \sigma_j (\beta)\vert} \right\rceil
= L+1 + \max_{2\leq j \leq d} \left\lceil \frac{\log n}{\log \vert \sigma_j (\beta)^{-1} \vert} \right\rceil
\end{equation*}
which we write as
\begin{equation*}
k = L+1 + \delta \log n + O(1), \quad \text{where} \quad \delta := \max_{2\leq j \leq d} \frac{1}{\log \vert \sigma_j(\beta)^{-1}\vert}.
\end{equation*}
Note that the $O(1)$ constant coming from the ceiling lies in $[0,1)$. For this choice of $k$,  $T_\beta^k(\frac{n}{\beta^{L+1}}) \in Y$.

Let $W$ be the maximum length of the $\beta$-expansions of the elements in $Y$. We therefore obtain asymptotic bounds of $R$,
\begin{equation*}
\delta \log n \leq R \leq W + 1 + \delta \log n \leq \delta' \log n,
\end{equation*}
where $\delta'>\delta$ can be chosen arbitrarily close.
\end{Proof}

In the course of the proof of Theorem \ref{ourThm} we deal with a problem caused by a constant coming from the length of the words we want to patch together. This issue can be circumvented by dividing the words into smaller subwords and glueing them back together afterwards. The following lemma ensures normality of the resulting word when patched back together.

\begin{Lem}\label{PatchingOfNormalWords}
Let $v=v_1 v_2 v_3 \ldots$ and $w = w_1 w_2 w_3 \ldots$ be $\mu$-normal words such that $\Vert v_i \Vert = \Vert w_i \Vert$, $\Vert v_i \Vert \rightarrow \infty$ and assume that the quantities
$\frac{\Vert v_{N+1}\Vert}{\Vert v_1 \Vert + \ldots + \Vert v_N \Vert}$ and $\frac{N}{\Vert v_1 \Vert + \ldots + \Vert v_N \Vert}$ tend to zero as $N$ tends to infinity. Then the word 
\begin{equation*}
u = v_1 w_1 v_2 w_2 v_3 w_3 \ldots 
\end{equation*}
is $\mu$-normal.
\end{Lem}

\begin{Proof}
We work with a fixed finite string $d=d_1 \ldots d_k$ of length $k$. Since $v$ is $\mu$-normal, we have
\begin{equation*}
N_d(v,n) = \sum_{i=1}^N N_d(v_i) + O(N) + O(\Vert v_{N+1} \Vert) \longrightarrow \mu(d) n
\end{equation*}
as $n \rightarrow \infty$, where $N$ is chosen such that $\Vert v_1 \Vert + \ldots + \Vert v_N \Vert \leq n < \Vert v_1 \Vert + \ldots + \Vert v_{N+1} \Vert$. The $O(N)$ contribution comes from possible occurrences in-between two words $v_i$ and $v_{i+1}$. Hence by the assumptions
\begin{equation*}
\sum_{i=1}^N N_d(v_i) \longrightarrow \mu(d)n
\end{equation*}
as $n \rightarrow \infty$. Then, by arguing similarly,
\begin{equation*}
N_d(u, n) = \sum_{i=1}^N N_d(v_i) + \sum_{i=1}^N N_d(w_i) + O(N) + O(\Vert w_{N+1} \Vert) = \mu(d) \frac{n}{2} + \mu(d) \frac{n}{2} + o(n).
\end{equation*}
Here, $N$ is chosen such that $\Vert v_1 \Vert + \Vert w_1 \Vert + \ldots + \Vert v_N \Vert + \Vert w_N \Vert \leq n < \Vert v_1 \Vert + \Vert w_1 \Vert + \ldots + \Vert v_{N+1} \Vert + \Vert w_{N+1} \Vert$. This shows the normality of $u$.
\end{Proof}

To verify the conditions of Lemma \ref{PatchingOfNormalWords} in our application, we need some basic number theoretic input.

\begin{Lem}
Denote by $p_N$ the $N$-th prime number. We have for $N \rightarrow \infty$
\begin{equation*}
(1) \quad \frac{\log p_{N+1}}{\log p_1 + \ldots + \log p_{N}} \longrightarrow 0, \quad \text{and} \quad 
(2) \quad \frac{N}{\log p_1 + \ldots + \log p_{N}} \longrightarrow 0.
\end{equation*}
\end{Lem}

\begin{Proof}
This is a consequence of the prime number theorem. We have 
\begin{equation*}
\sum_{i=1}^{p_N} \log p_i = \theta(p_N) \sim p_N \sim N\log N.
\end{equation*}
\end{Proof}

With these preliminaries we can prove our theorem.\\

\begin{Proofof}{Theorem \ref{ourThm}}
The polynomial $f(n) = a_gn^g + \ldots + a_1 n + a_0$ behaves asymptotically like $a_g n^g$. Hence for any $\epsilon>0$ and any $n$ large enough
\begin{equation*}
f(n) \leq (1+\epsilon) a_g n^g \quad \text{and} \quad f(n) \geq (1-\epsilon) a_g n^g.
\end{equation*}
Thus a consequence of Lemma \ref{boundsgeneralcase}, we have the upper bounds
\begin{align*}
R(f(n)) & \leq R((1+\epsilon)a_gn^g) \\
	& \leq \delta'( \log(1+\epsilon) + \log a_g + g \log n) \\
	& \leq C' \log n
\end{align*}
for some constant $C'>\delta' g$ arbitrarily close and $n$ large enough. Similarly we obtain lower bounds of the form
\begin{align*}
R(f(n)) \geq C \log n,
\end{align*}
where $C<\delta g$ is a positive constant and can be chosen arbitrarily close if $n$ is assumed to be large enough.

A direct consequence is an asymptotic upper bound for the total length of $f(n)$ when written in base $\beta$:
\begin{align*}
\Vert (f(n))_\beta \Vert &= L(f(n))+1+R(f(n))\\
	 &\leq L((1+\epsilon)a_g n^g) + 1 + R((1+\epsilon)a_g n^g) \\
	 &\leq \frac{g\log n + \log a_g}{\log \beta} + O(1) + C' \log n\\ 
	 &\leq C \frac{\log n}{\log \beta},
\end{align*}
for some (other) constant $C$. Note that $C$ only depends of $f$ and $\beta$.

However, applying Theorem \ref{MainThmABM} directly does not work as we do not have control of the size of the constant $C$. This can be avoided by choosing an integer $m\geq0$ such that $C/2^m \leq 1$ and dividing the words $(f(n))_\beta$ in $2^m$ words of (almost) equal length. Then we can apply Theorem \ref{MainThmABM} to show normality of the concatenations of those shorter words. They can subsequently be patched back together applying Lemma \ref{PatchingOfNormalWords} multiple times. Note that the lower bounds for $R$ enable us to use Lemma \ref{PatchingOfNormalWords}.

For a prime number $p$ we have $\Vert (f(p))_\beta \Vert \leq C \frac{\log p}{\log \beta}$, so $\Vert (f(p))_\beta \Vert \leq N$ is implied by $C \frac{\log p}{\log \beta} \leq N$. This is equivalent to
\begin{equation*}
p \leq \beta^{N/C}.
\end{equation*}
Thus, counting primes below $\beta^{N/C}$,
\begin{equation*}
\pi(\beta^{N/C}) \sim \frac{\beta^{N/C}}{N/C \log \beta} \geq \frac{C}{\log \beta} \beta^{N(C^{-1} - \epsilon)},
\end{equation*}
for any $\epsilon <1$ arbitrarily close and $N$ large enough. Here we see why we require $C\leq 1$, namely so that the condition of Theorem \ref{MainThmABM},
\begin{equation*}
\pi(\beta^{N/C}) \geq (\beta^N)^{1-\epsilon}
\end{equation*}
for any $\epsilon' > 0$, is implied by $\beta^{N(C^{-1} - \epsilon)}$ being eventually greater than $(\beta^N)^{1-\epsilon}$. Inserting the intermediary word to obtain admissibility in base $\beta$ does not destroy the normality of the sequence since we are inserting a word of constant length.
\end{Proofof}

\section{Special Case: Golden Ratio Base}

Let $\phi = \frac{1+\sqrt{5}}{2}$ be the golden ratio, i.e. the dominating root of the polynomial $x^2-x-1$. All positive integers have finite $\phi$-expansion (see for example Theorem 2 of \cite{frougny1992finite}). Considering this special case is interesting insofar that we can provide an \emph{exact} formula on the length of the fractional part. Let $n$ be a positive integer and denote by $L+1$ and $R$ the lengths of its integer and fractional part when written in base $\phi$. From the greedy algorithm we already know that
\begin{equation*}
\phi^L \leq n < \phi^{L+1} \Rightarrow L = \lfloor \frac{\log n}{\log \phi} \rfloor.
\end{equation*}
The following description of the fractional part is borrowed from \cite{GrabnerProdinger}.

\begin{Lem}\label{fracPart}
We have 
\begin{equation*}
R = \begin{cases} L & ,L \text{ even} \\ 
	L+1 & ,L \text{ odd.} \end{cases}
\end{equation*}
\end{Lem}

\begin{Proof}
Recall that $-\frac{1}{\phi}$ is the conjugate of $\phi$. Let $\sigma_2$ be the corresponding conjugation. Since it fixes the integers, we have
\begin{equation*}
n = \sigma_2(n) = \sum_{k= -R}^L \epsilon_k (-\frac{1}{\phi})^k = \sum_{k=-L}^R (-1)^k \epsilon_{-k} \phi^k.
\end{equation*}
The last sum can only be positive when $R$ is even. In this case we have $n < \phi^R + \phi^{R-2} + \phi^{R-4} + \ldots = \phi^{R+1}$ and since $\epsilon_{-(R-1)}=0$ also $n> \phi^{R} - \phi^{R-3} - \phi^{R-5} - \ldots = \phi^{R-1}$. This together with $\phi^L \leq n < \phi^{L+1}$  shows that $R$ is the even integer satisfying $L \leq R \leq L+1$.
\end{Proof}

The normality to base $\phi$ of the word
\begin{equation*}
(2)_\phi 0 (3)_\phi 0 (5)_\phi 0 \ldots
\end{equation*}
is then shown by first applying the generalized Copeland-Erd\H{o}s criterion Theorem \ref{MainThmABM} to the concatenations of only fractional or integer part and subsequently patching them back together using Lemma \ref{PatchingOfNormalWords} once.\\

\textbf{Acknowledgements:}
The author would like to thank Robert Tichy and Manfred Madritsch for their support, supervision and discussions about the subject. Special thanks to Tom\'{a}\v{s} V\'{a}vra who showed me \cite{GrabnerProdinger} and to Wolfgang Steiner who very patiently explained to me the proof of Proposition 3.5 in \cite{FrougnySteiner}. Part of this research was carried out at the Universit\'{e} de Lorraine in Nancy and at the Max-Planck Institute for Mathematics in Bonn and the author would like to thank these institutions for their hospitality.

\bibliographystyle{plain}
\bibliography{refs}

\end{document}